\documentclass{elsart}

\usepackage{amsfonts}		
\usepackage{latexsym}		

\let\Graphics0			
\usepackage{epsfig}
\newcommand{\XFig}[3]{\epsfxsize=#3\epsffile[#2]{#1}}
\newcommand{\Math}{{\slshape Mathematica}}

\newcommand{\RR}{{\mathbb R}}
\newcommand{\CC}{{\mathbb C}}
\newcommand{\HH}{{\mathbb H}}
\newcommand{\OO}{{\mathbb O}}

\newcommand{\nxn}{$n \times n$}
\newcommand{\Tr}{{\rm tr\,}}
\renewcommand{\bar}{\overline}
\renewcommand{\Re}{{\rm Re}}

\newcommand{\oi}{\,\hbox{\boldmath $i$}}
\newcommand{\oj}{\,\hbox{\boldmath $j$}}
\newcommand{\ok}{\,\hbox{\boldmath $k$}}
\newcommand{\ol}{\,\hbox{\boldmath $\ell$}}
\newcommand{\okl}{\,\hbox{\boldmath $k\ell$}}
\newcommand{\ojl}{\,\hbox{\boldmath $j\ell$}}
\newcommand{\oil}{\,\hbox{\boldmath $i\ell$}}

\newcommand\tttilde{\char'176}

\begin{document}


\begin{frontmatter}

\title{\bfseries Finding Octonionic Eigenvectors\\ Using \Math}

\author{Tevian Dray}
\address{
Department of Mathematics, Oregon State University,
                Corvallis, OR  97331, USA \\
{\tt tevian{\rm @}math.orst.edu} \\
}
\author{Corinne A. Manogue}
\address{
Department of Physics, Oregon State University,
                Corvallis, OR  97331, USA \\
{\tt corinne{\rm @}physics.orst.edu} \\
}

\begin{abstract}
The eigenvalue problem for $3\times3$ octonionic Hermitian matrices contains
some surprises, which we have reported elsewhere \cite{Eigen}.  In particular,
the eigenvalues need not be real, there are 6 rather than 3 real eigenvalues,
and the corresponding eigenvectors are not orthogonal in the usual sense.  The
nonassociativity of the octonions makes computations tricky, and all of these
results were first obtained via brute force (but exact) \Math\ computations.
Some of them, such as the computation of real eigenvalues, have subsequently
been implemented more elegantly; others have not.  We describe here the use of
\Math\ in analyzing this problem, and in particular its use in proving a
generalized orthogonality property for which no other proof is known.
\end{abstract}


\journal{Computer Physics Communications}
\date{17 June 1998}

\end{frontmatter}

\section{INTRODUCTION}

Finding the eigenvalues and eigenvectors of a given matrix is one of the basic
techniques in linear algebra, with countless applications.  The familiar case
of Hermitian (complex) matrices is very important, for instance in quantum
mechanics, where the fact that such matrices have real eigenvalues allows them
to represent physically observable quantities.

The eigenvalue problem is usually formulated over the complex numbers $\CC$,
including the reals $\RR$ as a special case.  In recent work \cite{Eigen}, we
considered the generalization to the other normed division algebras, namely
the quaternions $\HH$ and the octonions $\OO$.  Most of the basic properties
are retained, provided they are reinterpreted to take into account the lack of
commutativity of $\HH$ and $\OO$, and the lack of associativity of $\OO$.
However, there are a number of surprises, including the fact that such
matrices admit non-real eigenvalues.

Our most important results concern the eigenvalue problem for $3\times3$
octonionic Hermitian matrices, the {\em Jordan matrices\/}.  It turns out
\cite{Eigen} that such matrices admit more than the expected 3 real
eigenvalues, and that eigenvectors corresponding to different eigenvalues fail
in general to be orthogonal in the usual sense, although they do seem to be
orthogonal in a generalized sense.

Because of the lack of both commutativity and associativity, working with
octonionic matrices is rather tricky.  All of the above results were initially
discovered using a \Math\ package we have developed over the years for just
this purpose.  While we have subsequently been able to derive more elegant
derivations by hand for some of these results, we have not succeeded in doing
this for all.  In particular, the only current proof of the generalized
orthogonality property consists of a lengthy, brute force \Math\ computation,
which used 6 hours of CPU time on a Sparc20 with 224 Mb of memory.  This
article describes both the mathematics behind this result, and the \Math\
computation used to obtain it.

We set the stage in Section~\ref{Background} by reviewing some basic
properties about both the standard eigenvalue problem and the octonions, as
well as introducing our \Math\ package.  We then summarize the theory of
$3\times3$ octonionic Hermitian matrices in Section~\ref{JordanMatrices},
pointing out that the only known proof of some of the results in this section
involves direct computation using \Math.  In Section~\ref{Example} we give an
explicit example, and we discuss our results in Section~\ref{Discussion}.

\section{BACKGROUND}
\label{Background}

\subsection{The Standard Eigenvalue Problem}

We begin by collecting some standard results about the standard eigenvalue
problem.  We give the details of some of the proofs in order to emphasize the
use of both the commutativity and associativity of $\CC$.

The eigenvalue problem as usually stated is to find solutions $\lambda,v$ to
the equation
\begin{equation}
\label{Orig}
A v = \lambda v
\end{equation}
for a given square matrix $A$.  The basic properties of the eigenvalue problem
for \nxn\ complex Hermitian matrices are well-understood.

\begin{lem}
An \nxn\ complex Hermitian matrix $A$ has $n$ eigenvalues (counting
multiplicity), all of which are real.
\end{lem}

\begin{pf}
We give here only the proof that the eigenvalues are real.  Let $A$,
$v$, $\lambda$ satisfy (\ref{Orig}), with $A^\dagger=A$.  Then
\begin{equation}
\bar\lambda v^\dagger v
   = (A v)^\dagger v
   = v^\dagger A v
   = \lambda v^\dagger v
\end{equation}
so that if $v\ne0$ we have $v^\dagger v\ne0$, which forces
$\bar\lambda=\lambda$.
\qed
\end{pf}

\begin{lem}
Eigenvectors of an \nxn\ complex Hermitian matrix $A$ corresponding
to different eigenvalues are orthogonal.
\end{lem}

\begin{pf}
For $m=1,2$, let $v_m$ be an eigenvector of $A=A^\dagger$ with
eigenvalue $\lambda_m$.  By the previous lemma, $\lambda_m\in\RR$.  Then
\begin{equation}
\lambda_1 v_1^\dagger v_2
  = (A v_1)^\dagger v_2
  = v_1^\dagger A v_2
  = \lambda_2 v_1^\dagger v_2
\end{equation}
Then either $\lambda_1=\lambda_2$ or $v_1^\dagger v_2=0$.
\qed
\end{pf}

\begin{lem}
For any \nxn\ complex Hermitian matrix $A$, there exists an
orthonormal basis of $\CC^n$ consisting of eigenvectors of $A$.
\end{lem}

\begin{pf}
If all eigenvalues have multiplicity one, the result follows from the
previous lemma.  But the Gram-Schmidt orthogonalization process can be used on
any eigenspace corresponding to an eigenvalue with multiplicity greater than
one.
\qed
\end{pf}

\noindent
These lemmas are equivalent to the standard result that a complex Hermitian
matrix can always be diagonalized by a unitary transformation.  It is
important for what follows to realize that the form of the proofs given above
relies on both the commutativity and the associativity of $\CC$.

Combining the above results, it is easy to see that any (complex) Hermitian
matrix $A$ admits a decomposition in terms of an orthonormal basis
of eigenvectors.

\begin{thm}
Let $A$ be an \nxn\ complex Hermitian matrix.  Then $A$ can be
expanded as
\begin{equation}
A = \sum_{m=1}^n \lambda_m v_m v_m^\dagger
\label{Decomp}
\end{equation}
where $\{v_m; ~m=1,...,n\}$ is an orthonormal basis of eigenvectors
corresponding to eigenvalues $\lambda_m$.
\end{thm}

\begin{pf}
By the previous lemma, there exists an orthonormal basis $\{v_m\}$ of
eigenvectors.  It then suffices to check that
\begin{equation}
\sum_{m=1}^n \lambda_m v_m v_m^\dagger v_k = \lambda_k v_k
\end{equation}
But this follows by direct computation using orthonormality.
\qed
\end{pf}

\noindent
Furthermore, the set of eigenvalues $\{\lambda_m\}$ is unique, and the (unit)
eigenvectors are unique up to unitary transformations in the separate
eigenspaces (which reduce to multiplication by a complex phase for eigenvalues
of multiplicity one).

\subsection{Octonions}

The quaternions $\HH$ double the dimension of the complex numbers by adding
two additional square roots of $-1$, usually denoted $j$ and $k$.  The
multiplication table follows from
\begin{equation}
i^2=j^2=k^2=-1 \qquad ij=k=-ji
\end{equation}
and associativity; note that $\HH$ is not commutative.  Equivalently, $\HH$
can be viewed as the sum of 2 copies of the complex numbers
\begin{equation}
\HH = \CC + k\CC
\end{equation}
with $j$ being defined by $j=ki$.

\begin{figure}
\begin{center}
\XFig{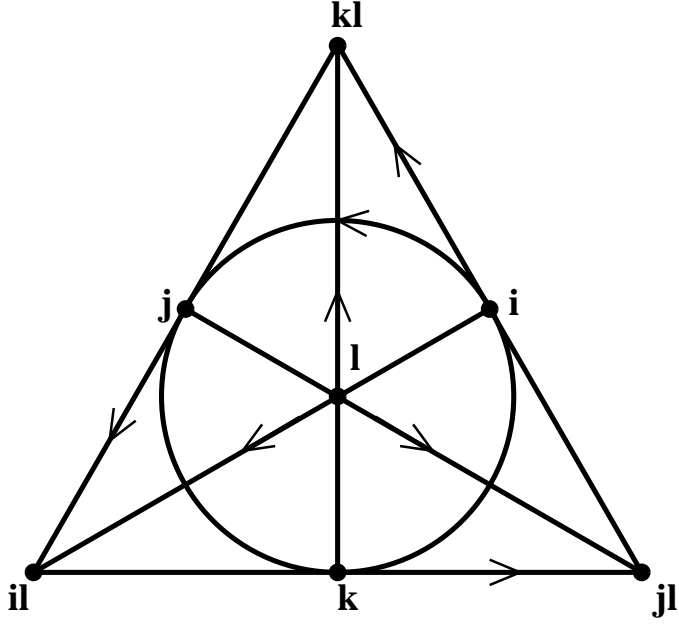}{68 168 543 614}{3.5in}
\end{center}
\caption{The representation of the octonionic multiplication table using the
7-point projective plane.  Each of the 7 oriented lines gives a quaternionic
triple.}
\label{Omult}
\end{figure}

The octonions $\OO$ in turn can be viewed as the direct sum of two copies of
the quaternions
\footnote{This construction of a new division algebra from 2 copies of another
is a special case of the Cayley-Dickson process; for modern treatments, see
\cite{ReeseHarvey,Kantor,LounestoI,Wene,Schafer}.}
\begin{equation}
\OO = \HH + \HH\ell = (\CC + k\CC) + (\CC + k\CC)\ell 
\end{equation}
where $\ell$ is yet another square root of $-1$.  The octonions are thus
spanned by the identity element $1$ and the 7 imaginary units
$\{i,j,k,k\ell,j\ell,i\ell,\ell\}$.  These units can be grouped into (the
imaginary parts of) quaternionic subspaces in 7 different ways; these will be
referred to as ``triples''.  Any three of these imaginary units which do not
lie in a such a triple anti-associate.  The multiplication table can be neatly
summarized by appropriately labeling the 7-point projective plane, as shown in
Figure~\ref{Omult}.

Even though the octonions are not associative, since any 2 octonions lie in a
quaternionic subspace, products involving only 2 different octonions (and
their octonionic conjugates) do associate.  For example,
\begin{equation}
p(pq) = p^2 q
\end{equation}
which is a weak form of associativity known as {\em alternativity\/}.

We use the notation $\bar{a}$ to denote the (octonionic) conjugate of the
octonion $a$,
\begin{equation}
|a|^2 := a\bar{a}
\end{equation}
to denote the squared norm of $a$, $A^\dagger$ to denote the (octonionic)
Hermitian conjugate of the matrix $A$,
\begin{equation}
[a,b] := ab - ba
\end{equation}
to denote the commutator of $a$ and $b$, and
\begin{equation}
[a,b,c] := (ab)c - a(bc)
\end{equation}
to denote the associator of $a$, $b$, $c$.  Both the commutator and the
associator are purely imaginary, totally antisymmetric, and change sign if any
one of their arguments is replaced by its conjugate.  Another octonionic
product with these properties is given by the {\em associative 3-form\/}
\cite{ReeseHarvey,Gureirch}
\begin{equation}
\Phi(a,b,c)
  = {1\over2} \Re \Big( a(\bar{b}c) - c(\bar{b}a) \Big)
\end{equation}
which reduces to the vector triple product when $a$, $b$, $c$ are imaginary
quaternions.
\begin{equation}
\Phi(a,b,c) = {1\over2} \, \Re( [a,\bar{b}] c)
\label{PhiEq}
\end{equation}

\subsection{\Math\ Package}
\label{Package}

We needed a way to easily manipulate octonions and octonionic matrices --- it
is quite difficult to unlearn associativity!  There are 2 complementary
approaches, depending on whether it is desired to manipulate abstract
octonions or whether an explicit basis can be used.  For our purposes, it was
initially quite sufficient to work with an abstract basis: We define an
octonion to be a list with 8 elements
\begin{eqnarray}
a &=& a_1 +a_2 \oi +a_3 \oj +a_4 \ok 
	-a_5 \okl -a_6 \ojl -a_7 \oil
	+a_8 \ol \nonumber\\
  &=& \{a_1,a_2,a_3,a_4,a_5,a_6,a_7,a_8\}
\end{eqnarray}
where the signs are conventional, and where we have emphasized the role of the
imaginary units as ``basis vectors'' by writing them in boldface.  Octonionic
multiplication can then be expressed as a series of rules of the form
\begin{equation}
\oi * \oj = \ok
\end{equation}
where we have introduced the symbol $*$ to represent octonionic
multiplication.  Octonionic conjugation can be defined as a series of rules of
the form
\begin{equation}
\bar{\hbox{\boldmath $i$}} = -\oi
\end{equation}
and everything else can be defined in terms of these 2 basic operations.

\begin{figure}
\begin{center}
\XFig{Find2.ps}{55 538 231 728}{2.5in}
\end{center}
\if\Graphics1\vspace{-0.4in}\fi
\caption{A sample \Math\ session illustrating basic manipulation of octonions
using an explicit basis representation, where $O1$, $O2$, $O3$ are defined as
in (\ref{Generic}).}
\label{BasicsII}
\vspace{0.2in}
\end{figure}


Furthermore, as already stated, computations involving small numbers of
octonions can be dramatically simplified.  For instance, 3 arbitrary octonions
can be assumed to take the form
\begin{eqnarray}
O1 &=& a_1 + a_2 \oi \nonumber\\
O2 &=& b_1 + b_2 \oi + b_3 \oj
\label{Generic}\\
O3 &=& c_1 + c_2 \oi + c_3 \oj + c_4 \ok + c_8 \nonumber\ol
\end{eqnarray}

When implementing these ideas using \Math, it turned out to be more efficient
to define the 2 fundamental operations directly on lists, rather than building
them up in terms of rules.  For instance, conjugation is more easily defined
by
\begin{eqnarray}
\tt bar[\{x1\_,x2\_,x3\_,x4\_,x5\_,x6\_,x7\_,x8\_\}]:= \nonumber\\
	\qquad\qquad\tt\{x1,-x2,-x3,-x4,-x5,-x6,-x7,-x8\}
\end{eqnarray}
and there is an analogous definition of octonionic multiplication, called {\tt
Omult}.  A illustration of computation in an explicit basis appears in
Figure~\ref{BasicsII}.

Matrices can now be constructed as nested lists, and operations such as matrix
multiplication can easily be defined.  However, matrix expressions involving
several octonions quickly become unwieldy.  We therefore supplemented the
above basic definitions in terms of an explicit basis with an alternative set
of definitions using abstract octonions.  In the process, we took full
advantage of the formatting capabilities in \Math\ 3.0.  An illustration of
computation using the abstract approach appears in Figure~\ref{BasicsI}.

\begin{figure}
\begin{center}
\XFig{Find1.ps}{55 469 331 728}{4in}
\end{center}
\if\Graphics1\vspace{-0.6in}\else\vspace{-0.4in}\fi
\caption{A sample \Math\ session illustrating basic manipulation of
abstract octonions.}
\label{BasicsI}
\end{figure}


Even though we chose not to implement the concrete version of the basic
octonionic operations in terms of rules, \Math's ability to manipulate such
rules was crucial in constructing the abstract version.  Especially when
teaching \Math\ to both manipulate and print matrices of lists (i.e.\
octonions), the ability to easily modify the code to recognize special cases
proved extremely helpful.
(Octonions are closely related to {\it Clifford algebras}; see for instance
\cite{Jorg}, in which the representation theory of Clifford algebras is
extended so as to employ octonions.  There are a number of packages for
manipulating Clifford algebras, some of which are described in
\cite{CliffordBook}.  Several of these can handle octonions, such as the
program CLICAL by Pertti Lounesto \cite{CLICAL} and the Maple package CLIFFORD
by Rafa\l \ Ab\l amowicz \cite{CLIFFORD}, both of which introduce
the octonions as paravectors over an appropriate Clifford algebra
\cite{Lounesto2}.  There is also an older Maple package Octonion (and the
companion package Clifford) by J\"org Schray \cite{Schray}.)

\newpage

\section{\boldmath $3\times3$ OCTONIONIC HERMITIAN MATRICES}
\label{JordanMatrices}

It is not immediately obvious that $3\times3$ octonionic Hermitian matrices
have a well-defined determinant, let alone a characteristic equation.  We
therefore first discuss some of the properties of these matrices before
turning to the eigenvalue problem.


\subsection{Jordan Matrices}

The $3\times3$ octonionic Hermitian matrices, henceforth referred to as the
{\em Jordan matrices\/}, form the {\em exceptional Jordan algebra\/} under
the {\em Jordan product\/}
\begin{equation}
A \circ B := {1\over2} (AB + BA)
\end{equation}
which is commutative, but not associative.  A special case of this is
\begin{equation}
A^2 \equiv A \circ A
\end{equation}
and we {\em define\/}
\begin{equation}
A^3 := A^2 \circ A = A \circ A^2
\end{equation}

Remarkably, with these definitions, Jordan matrices satisfy the usual
characteristic equation (see e.g.\ \cite{ReeseHarvey})
\begin{equation}
A^3 - (\Tr A) \, A^2 + \sigma(A) \, A - (\det A) \, I = 0
\label{Char}
\end{equation}
where $\sigma(A)$ is defined by
\begin{equation}
\sigma(A) := {1\over2} \left( (\Tr A)^2 - \Tr (A^2) \right)
\end{equation}
and where the determinant $\det A$ of $A$ is defined abstractly in terms of
the Freudenthal product \cite{Jordan,Freudenthal}
\begin{eqnarray}
\nonumber
A*B = A \circ B - {1\over2} \Big(A\,\Tr(B)+B\,\Tr(A)\Big)
                + {1\over2} \Big(\Tr(A)\,\Tr(B)-\Tr(A\circ B)\Big) \\
\end{eqnarray}
which leads to
\begin{equation}
\det(A) = {1\over3} \, \Tr \Big( (A*A) \circ A \Big)
\end{equation}
\goodbreak\noindent
Concretely, if
\begin{equation}
A = \pmatrix{p& a& \bar{b}\cr \bar{a}& m& c\cr b& \bar{c}& n\cr}
\label{Three}
\end{equation}
with $p,m,n\in\RR$ and $a,b,c\in\OO$, then
\begin{eqnarray}
\Tr A &=& p + m + n \\
\sigma(A) &=& pm + pn + mn - |a|^2 - |b|^2 - |c|^2 \\
\det A &=& pmn + b(ac) + \bar{b(ac)} - n|a|^2 - m|b|^2 - p|c|^2
\label{ThreeEq}
\end{eqnarray}
The determinant calculation is illustrated in Figure~\ref{MathMatrices}.

\begin{figure}
\begin{center}
\if\Graphics1
	\XFig{Find4.ps}{55 448 448 728}{5.55in}
\else
	\XFig{Find4.ps}{55 448 448 728}{5.4in}
\fi
\end{center}
\if\Graphics1\vspace{-0.3in}\fi
\caption{A \Math\ session showing the calculation of the determinant for a
generic Jordan matrix, both abstractly and in terms of an explicit basis.}
\label{MathMatrices}
\vspace{0.2in}
\end{figure}

\subsection{The Real Eigenvalue Problem}

Each division algebra can be rewritten as a real matrix algebra of appropriate
dimension (see e.g.\ \cite{Eigen}).  Under this identification, a Hermitian
matrix over any of the division algebras becomes a real symmetric matrix.  It
is therefore clear that a $3\times3$ octonionic Hermitian matrix must have
$8\times3=24$ real eigenvalues~\cite{Horwitz}.  However, as we now discuss,
instead of having (a maximum of) 3 distinct real eigenvalues, each with
multiplicity 8, there appear to be (a maximum of) 6 distinct real eigenvalues,
each with multiplicity 4.

The reason for this is that, somewhat surprisingly, a (real) eigenvalue
$\lambda$ of a Jordan matrix $A$ does {\em not\/} in general satisfy the
characteristic equation.  To see this, consider the eigenvalue equation
(\ref{Orig}), with $A$ as in (\ref{Three}), $\lambda\in\RR$, and where
$v=\pmatrix{x\cr y\cr z\cr}$.  Assuming without loss of generality that
$z\ne0$, explicit computation yields \cite{Eigen}
\begin{eqnarray}
\Big[ \det(\lambda I - A) \Big] z 
  &\equiv& \left[ 
        \lambda^3 - (\Tr A) \, \lambda^2 + \sigma(A) \, \lambda - \det A
     \right] z \nonumber\\
  &=& b \Big( a (cz) \Big) + \bar{c} \left( \bar{a} (\bar{b}z) \right)
     - \left[ b(ac) + (\bar{c}\,\bar{a})\bar{b} \right] z
\label{CharLam}
\end{eqnarray}

If $a$, $b$, $c$, and $z$ associate, the RHS of (\ref{CharLam}) vanishes, and
$\lambda$ does indeed satisfy the characteristic equation (\ref{Char}); this
will not happen in general.  However, since the LHS of (\ref{CharLam}) is a
real multiple of $z$, this must also be true of the RHS, so that
\begin{equation}
b \Big( a (cz) \Big) + \bar{c} \left( \bar{a} (\bar{b}z) \right)
     - \left[ b(ac) + (\bar{c}\,\bar{a})\bar{b} \right] z
  = rz \qquad\qquad r\in\RR
\label{RReal}
\end{equation}
which can be solved to yield a quadratic equation for $r$ as well as
constraints on $z$.

\begin{thm}[Dray \& Manogue \cite{Eigen}]
The real eigenvalues of the $3\times3$ octonionic Hermitian matrix
$A$ satisfy the modified characteristic equation
\begin{equation}
\det(\lambda I - A)
  = \lambda^3 - (\Tr A) \, \lambda^2 + \sigma(A) \, \lambda - \det A
  = r
\label{Lameq}
\end{equation}
where $r$ is either of the two roots of
\begin{equation}
r^2 + 4\Phi(a,b,c) \, r - \Big| [a,b,c] \Big|^2 = 0
\label{Req}
\end{equation}
with $a,b,c$ as defined by (\ref{Three}) and where $\Phi$ was defined in
(\ref{PhiEq}).
\label{CharThm}
\end{thm}

\noindent
Furthermore, provided that $[a,b,c]\ne0$, each of $x$, $y$, and $z$ can be
shown to admit an expansion of the form (given for $z$ only)
\begin{cor}
With $A$ and $r$ as above, and assuming $[a,b,c]\ne0$,
\begin{equation}
z = (\alpha a + \beta b + \gamma c + \delta)
      \left( 1 + {[a,b,c] \, r \over \Big| [a,b,c] \Big|^2} \right)
\label{Zeq}
\end{equation}
with $\alpha,\beta,\gamma,\delta\in\RR$.
\end{cor}

\begin{pf}
Both of these results were obtained using \Math\ to solve (\ref{RReal}) by
brute force for real $r$ and octonionic $z$ given generic octonions $a$, $b$,
$c$.  An outline of the computation appears in Figure~\ref{MathREq}.
\qed
\end{pf}

\begin{figure}
\begin{center}
\if\Graphics1
	\XFig{Find3.ps}{55 252 508 719}{6in}
\else
	\XFig{Find3.ps}{55 252 508 719}{5.45in}
\fi
\end{center}
\if\Graphics1\vspace{-0.55in}\else\vspace{-0.3in}\fi
\caption{An outline of the derivation of the modified characteristic equation
using \Math.  We first construct $Eq$ to be (\ref{RReal}), then solve half
these equations for some of the components of $z$ in terms of the others (this
is {\it ZRule}).  We then solve the remaining equations for $r$ (this is {\it
RRule}) and verify that the solutions have the form claimed in
Theorem~\ref{CharThm}, i.e.\ that they satisfy (\ref{Req}) (which is $Req$).
Finally, we verify that $z$ has the form claimed in
(\ref{Zeq}) (this is {\it PRule}).}
\label{MathREq}
\if\Graphics1\vspace{0.2in}\else\vspace{0.1in}\fi
\end{figure}

\noindent
The real parameters $\alpha$,$\beta$,$\gamma$,$\delta$ may be freely specified
for one (nonzero) component, say $z$; the remaining components $x$,$y$ have a
similar form.

The solutions of (\ref{RReal}) are real, since the corresponding $24\times24$
real symmetric matrix has 24 real eigenvalues.  We will refer to the 3 real
solutions of (\ref{RReal}) corresponding to a single value of $r$ as a {\em
family\/} of eigenvalues of $A$.  There are thus 2 families of real
eigenvalues, each corresponding to a 4 independent (over $\RR$) eigenvectors.

We note several intriguing properties of these results.  If $A$ is in fact
complex, then the only solution of (\ref{Req}) is $r=0$, and we recover the
usual characteristic equation with a unique set of 3 (real) eigenvalues.  If
$A$ is quaternionic, then one solution of (\ref{Req}) is $r=0$, leading to the
standard set of 3 real eigenvalues and their corresponding quaternionic
eigenvectors.  However, unless $a$, $b$, $c$ involve only two independent
imaginary quaternionic directions (in which case $\Phi(a,b,c)=0=[a,b,c]$),
there will also be a nonzero solution for $r$, leading to a second set of 3
real eigenvalues.  Finally, if $A$ is octonionic (so that in particular
$[a,b,c]\ne0$), then there are two distinct solutions for $r$, and hence two
different sets of real eigenvalues, with corresponding eigenvectors.  Note
that if $\det{A}=0\ne[a,b,c]$ then all of the eigenvalues of $A$ will be
nonzero!

\subsection{Orthogonality}

The final surprise lies with the orthogonality condition for eigenvectors $v,w$
corresponding to different eigenvalues.  It is {\em not\/} true that
$v^\dagger w=0$, although the real part of this expression does vanish
\cite{Eigen}.  But, at least in the $2\times2$ case \cite{Eigen}, it is
straightforward to show that what is needed to ensure a decomposition of the
form (\ref{Decomp}) is the following generalized notion of orthogonality,
which does in fact hold.

\begin{defn}[Dray \& Manogue \cite{Eigen}]
Let $v$ and $w$ be two octonionic vectors.  We will say that $w$ is
{\em orthogonal} to $v$ if
\begin{equation}
\label{Ortho}
(v v^\dagger) \, w = 0
\end{equation}
\end{defn}

In the $3\times3$ case, a lengthy, direct computation verifies that
eigenvectors with different real eigenvalues satisfy (\ref{Ortho}) {\em
provided that\/} the same value of $r$ is used for both eigenvectors.

\begin{thm}[Dray \& Manogue \cite{Eigen}]
If $v$ and $w$ are eigenvectors of the $3\times3$ octonionic Hermitian matrix
$A$ corresponding to different real eigenvalues in the same family (same $r$
value), then $v$ and $w$ are mutually orthogonal in the sense of
(\ref{Ortho}).
\label{OrthoThm}
\end{thm}

\begin{pf}
The modified characteristic equation (\ref{Lameq}) can be used to
eliminate cubic and higher powers of $\lambda$ from any expression.
Furthermore, given two distinct eigenvalues $\lambda_1\ne\lambda_2$,
subtracting the two versions of (\ref{Lameq}) and factoring the result leads
to the equation
\begin{equation}
(\lambda_1^2+\lambda_1\lambda_2+\lambda_2^2)
  - (\Tr{A}) (\lambda_1+\lambda_2) + \sigma(A)
  = 0
\label{TwoLamEq}
\end{equation}
which can be used to eliminate quadratic terms in one of the eigenvalues.  We
used \Math\ to implement these simplifications in a brute force verification
of (\ref{Ortho}) in this context, which used 6 hours of CPU time on a SUN
Sparc20 with 224 Mb of RAM.  A summary of the computation appears in 
Figures \ref{MathOrthoI} and \ref{MathOrthoII}.
\qed
\end{pf}

\noindent
For Jordan matrices, we thus obtain {\em two\/} decompositions of the form
(\ref{Decomp}), corresponding to the two sets of real eigenvalues.  For each,
the eigenvectors are fixed up to orthogonal transformations which preserve the
form (\ref{Zeq}) of $z$.

\begin{figure}
\begin{center}
\XFig{FindOrtho1.ps}{76 214 500 716}{5in}
\end{center}
\if\Graphics1\vspace{-0.3in}\fi
\caption{The construction of eigenvectors $V[i]$ of $A$ with different real
eigenvalues $q[i]$.  The characteristic equation (\ref{Lameq}) is implemented
via {\it q3Rule}; the condition (\ref{PhiEq}) on $r$ is given by {\it r2Rule},
and {\it qRule} is the extra condition (\ref{TwoLamEq}).  Finally,
$VV=V[1]\,V[1]^\dagger$.}
\label{MathOrthoI}
\if\Graphics1\vspace{0.5in}\else\vspace{0.2in}\fi
\end{figure}

\begin{figure}
\begin{center}
\XFig{FindOrtho2.ps}{76 187 500 716}{5in}
\end{center}
\if\Graphics1\vspace{-0.5in}\fi
\caption{A summary of the \Math\ 2.2 computation used to show that
eigenvectors $V$, $W$ of Jordan matrices with different real eigenvalues are
orthogonal in the generalized sense.  Note that each term was simplified
separately.}
\label{MathOrthoII}
\vspace{1.5in}
\end{figure}

\begin{thm}[Dray \& Manogue \cite{Eigen}]
Let $A$ be a $3\times3$ octonionic Hermitian matrix.  Then $A$ can be expanded
as in (\ref{Decomp}), where $\{v_1, v_2, v_3\}$ is an orthonormal basis, as per
(\ref{Ortho}), of eigenvectors of $A$ corresponding to the real eigenvalues
$\lambda_m$, which belong to the same family (same $r$ value).
\label{DecompThm}
\end{thm}

\begin{pf}
Fix a family of real eigenvalues of $A$ by fixing $r$.  If the eigenvalues are
distinct, then the previous theorem guarantees the existence of an orthonormal
basis of eigenvectors, and the result follows.  \hfill\break\indent If the
eigenvalues are the same, the family consists of a single real eigenvalue
$\lambda$ with multiplicity 3.  Then $\Tr(A)=3\lambda$ and
$\sigma(A)=3\lambda^2$.  Writing out these two equations in terms of the
components (\ref{Three}) of $A$, and inserting the first into the second,
results in a quadratic equation for $\lambda$; the discriminant $D$ of this
equation satisfies $D\le0$.  But $\lambda$ is assumed to be real, which forces
$D=0$, which in turn forces $A$ to be a multiple of the identity matrix, for
which the result holds.  \hfill\break\indent The remaining case is when one
eigenvalue, say $\mu$, has multiplicity 2 and one has multiplicity 1.  Letting
$v$ be a (normalized) eigenvector with eigenvalue $\mu$, consider the matrix
\begin{equation}
X = A - \alpha \> v v^\dagger
\label{Two}
\end{equation}
with $\alpha\in\RR$.  For most values of $\alpha$, $X$ will have 3 distinct
real eigenvalues, whose eigenvectors will be orthogonal by the previous
theorem.  But this means that eigenvectors of $X$ are also eigenvectors of
$A$; the required decomposition of $A$ is obtained from that of $X$ simply by
solving (\ref{Two}) for $A$.
\qed
\end{pf}

\noindent
Note in particular that for some quaternionic matrices with determinant equal
to zero, one and only one of these two decompositions will contain the
eigenvalue zero.

\section{Example}
\label{Example}

Let $s$ be given by 
\begin{equation}
s=\cos\theta+k\ell\,\sin\theta
\end{equation}
and consider the matrix
\begin{equation}
B = \pmatrix{~~p & ~~iq & kqs \cr
                -iq & ~~p & jq\cr
                -kqs & -jq & p \cr}
\label{BDef}
\end{equation}
noting that $B$ is quaternionic if $\theta=0$.  Turning first to the
equation for $r$, (\ref{RReal}) becomes
\begin{equation}
r^2 + 4 q^3 r \cos\theta - 4 q^6 \sin^2\theta = 0
\end{equation}
with solutions
\begin{equation}
r_\pm = - 2 q^3 ( \cos\theta \pm 1 )
\end{equation}
Since
\begin{eqnarray}
\Tr B &=& 3p \\
\sigma(B) &=& 3 (p^2-q^2) \\
\det B &=& p^3 - 3 p q^2 + 2 q^3 \cos\theta
\end{eqnarray}
the eigenvalue equation (\ref{Lameq}) becomes 
\begin{eqnarray}
0
  &=& \lambda^3 - 3p \, \lambda^2
        + 3(p^2-q^2) \, \lambda - (p^3 - 3 p q^2 \mp 2 q^3) \\
  &=& (\lambda - p \mp q)^2 \, (\lambda - p \pm 2q)
\end{eqnarray}

\begin{figure}
\null\vspace{-5pt}
\begin{center}
\XFig{FindEx.ps}{55 125 416 728}{4.75in}
\end{center}
\if\Graphics1\vspace{-0.4in}\else\vspace{-0.2in}\fi
\caption{A \Math\ computation illustrating that the vectors $u_+$, $v_+$,
$w_+$ given in (\ref{EigenI})--(\ref{EigenIII}) are indeed eigenvectors of the
Jordan matrix $B$ given in (\ref{BDef}) with the given eigenvalues (with
$f_u=1$, $f_v=\oi$, and $f_w=\oj$), that these eigenvectors are only
orthogonal in the generalized sense of Theorem~\ref{OrthoThm}, and that
they lead to a decomposition of $B$ as implied by Theorem~\ref{DecompThm}.
(Normalization factors have been added in the final computation, since the
given vectors are not unit vectors.)}
\label{MathEx}
\end{figure}

An orthonormal basis of eigenvectors associated with these eigenvalues is
\begin{eqnarray}
\label{EigenI}
\lambda_u &= p \pm q: \quad
        u_\pm &= \pmatrix{i\cr 0\cr j\cr}
                (f_u R_\pm) \\
\noalign{\smallskip}
\lambda_v &= p \pm q: \quad
        v_\pm &= \pmatrix{j\cr 2ks\cr i\cr}
                (f_v R_\pm) \\
\noalign{\smallskip}
\label{EigenIII}
\lambda_w &= p \mp 2q: \quad
        w_\pm &= \pmatrix{~j\cr -ks\cr ~i\cr}
                (f_w R_\pm)
\end{eqnarray}
where $f_u$, $f_v$, $f_w$ are arbitrary linear
combinations of $a=iq$, $b=-kqs$, $c=jq$, i.e.\
\begin{equation}
  f_u,f_v,f_w \in \langle 1,a,b,c \rangle
\end{equation}
and where $R_\pm$ is given by
\begin{equation}
R_\pm = \cases
  {\sin{\theta\over2}+k\ell\,\cos{\theta\over2}\cr
   \noalign{\smallskip}
   \cos{\theta\over2}-k\ell\,\sin{\theta\over2}\cr
  }
\end{equation}
Note that in the limiting case $\theta\to0$, each $f$ is quaternionic, and $R_-$
reduces to $1$ while $R_+$ becomes $k\ell$.

As expected, {\em provided\/} one fixes a family of eigenvectors and
eigenvalues arising from a given choice of $r$, these eigenvectors satisfy the
orthogonality property (\ref{Ortho}) and thus lead to a decomposition of the
form (\ref{Decomp}).  A partial verification of this using \Math\ is given in
Figure~\ref{MathEx}.

\section{Discussion}
\label{Discussion}

There are 2 quite different surprising aspects of our work: the mathematical
changes needed to extend the eigenvalue problem to the octonions, and the fact
that we have only been able to prove one of our key results using computer
algebra.  We discuss each of these in turn.

It is of course intriguing that the eigenvalue problem over the octonions
changes so much, for instance in that there are unexpectedly many real
eigenvalues.  But we find it remarkable that so much of the standard structure
remains, provided it is reinterpreted appropriately.  The most striking
example of this is the need to generalize what is meant by orthogonality.

We can relate our notion of orthonormality to the usual one by noting that a
basis of $\OO^n$ which is orthonormal in the sense (\ref{Ortho}) satisfies
\begin{equation}
vv^\dagger + ... + ww^\dagger=I
\end{equation}
which follows directly from the definition.  If we define a matrix $Q$ whose
columns are just $v,...,w$, then this statement is equivalent to
\begin{equation}
Q Q^\dagger = I
\end{equation}
Over the quaternions, left matrix inverses are the same as right matrix
inverses, and we would also have
\begin{equation}
Q^\dagger Q = I
\end{equation}
or equivalently
\begin{equation}
v^\dagger v = 1 = ... = w^\dagger w; \qquad v^\dagger w = 0 = ...
\end{equation}
which is just the standard notion of orthogonality.  These two notions of
orthogonality fail to be equivalent over the octonions; we have been led to
view the former as more fundamental.

Turning to our proof-by-computer, we reiterate that the only proof we
currently have of our main orthogonality result, namely
Theorem~\ref{OrthoThm}, uses \Math\ to explicitly perform a horrendous, but
exact, algebraic computation.  While one could hope for a more elegant
mathematical proof of this result, the \Math\ computation nevertheless
establishes a result which would otherwise remain for the moment merely a
conjecture.  This is a good example of being able to use the computer to
verify one's intuition when it may not be possible to do so otherwise.

But even more is true: Throughout our work with the octonions, the ability to
manipulate octonionic expressions quickly and accurately has been crucial in
{\em developing} our intuition in the first place.  We do not feel that we
would have been able to reach anything like our current understanding of the
applications of the octonions to physics, on which we continue to be working
actively, without the availability of a package such as the one described
here.

Finally, this computation was initially done in 1996 using \Math\ 2.2.  While
preparing this paper, we attempted to reproduce the computation using \Math\
3.0 --- and couldn't!  Even for \Math\ 2.2, it was necessary to massage the
computation by hand in order to succeed.  One way this was done (between
Figure~\ref{MathOrthoI} and Figure~\ref{MathOrthoII}) was by saving some
intermediate steps to files and then restarting the kernel.  Another technique
was not to simplify all the components of an expression at the same time.  For
instance, in Figure~\ref{MathOrthoII}, the vector ${\tt VVW}=(V V^\dagger)W$
contains 3 octonions, each of which requires roughly 8 Mb.

Comparing the computations both versions of \Math\ could handle, \Math\ 3.0
appears to require nearly 4 times as much CPU time for the same computation;
this was for identical inputs, with a minimum of special formatting.  It is
unfortunate that the many nice features of \Math\ 3.0 appear to require such a
high price.


\end{document}